\input harvmac
\font\tit=cmr12
\font\male=cmr8
\hfill MPI-PhT/95-77
\vskip 3truecm
\centerline{\tit QUANTUM $E(2)$ GROUPS}
\bigskip
\centerline{\tit AND LIE BIALGEBRA STRUCTURES}
\vskip 2truecm
\centerline{J. Sobczyk\footnote*{{\male On leave on absence from
Institute for Theoretical Physics, Wroc\l aw University, Poland, E-mail:
jsobczyk@proton.ift.uni.wroc.pl}}}
\smallskip
\centerline{Max-Planck-Institut f\"ur Physik, Werner Heisenberg
Institute}
\centerline{F\"ohringer Ring 6, D-80805 Munich, Germany}
\vskip 2truecm
\centerline{ABSTRACT}
\smallskip
Lie bialgebra structures on
$e(2)$ are classified. For two Lie bialgebra structures which are not
coboundaries
(i.e. which are not determined by a classical $r$-matrix) we solve the
cocycle condition, find the Lie-Poisson brackets and obtain quantum group
relations.
There is one to one correspondence between Lie bialgebra structures on
$e(2)$
and possible quantum deformations of $U(e(2))$ and $E(2)$.
\bigskip\noindent
revised

\noindent
October 1995\hfill\break
\pageno=0
\vfill\eject
{\bf 1.} Quantum deformations
\ref\drber{V.G. Drinfeld, {\it Quantum groups}, Proc. ICM Berkeley 1986,
vol 1, 789}
\ref\frt{L.D. Faddeev, N.Yu. Reshetikhin and L.A. Takhtajan, Leningrad
Math. J. {\bf 1} (1990) 193}
of the $D=2$ Euclidean group $E(2)$ and its
universal enveloping algebra $U(e(2))$ turn out to be
a useful laboratory to study various aspects of quantum groups
\ref\cgsta{E. Celeghini, R. Giachetti, R. Sorace and M. Tarlini, J.
Math. Phys. {\bf 31} (1990) 2548} - \nref\wor
{S.L. Woronowicz, Lett. Math. Phys. {\bf 23} (1991) 251; S.L.
Woronowicz, Comm. Math. Phys. {\bf 149} (1992) 637}
\ref\swz{P. Schupp, P. Watts and B. Zumino, Lett. Math. Phys. {\bf 24}
(1992) 141}.
It is
one of the simplest examples of non-simple Lie group and there is no
canonical way to introduce its deformation. During last five
years many approaches have been developed
\ref\masa{P. Ma\'{s}lanka, J. Math. Phys. {\bf 35} (1994) 1976} -
\ref\bcgst{A. Ballesteros, E. Celeghini, R. Giachetti, R. Sorace and M.
Tarlini, J. Phys. A {\bf 26} (1993) 7495}
to construct such deformations.
The study of $E_q(2)$ is interesting for itself
(one
can ask questions about how many different quantum deformations exist in
this case, about classical $r$ and quantum $R$ matrices,
differential calculi, representations, bicrossproduct structures etc.)
but is also useful in order to understand properties of quantum deformations
of $D=4$ Poincar\'{e} group
\ref\lnrt{J. Lukierski, A. Nowicki, H. Ruegg and V. Tolstoy, Phys. Lett.
B {\bf 264} (1991) 331}. The structure of $D=4$ quantum groups
are important to explore if one wishes
to examine possible implications of quantum groups ideas in physics.
Recently interesting results were obtained in this direction including a
classification of possible deformations of $D=4$ Poincar\'{e} group
\ref\pw{P. Podle\'{s} and S.L. Woronowicz, {\it On the classification of
quantum Poincar\'{e} groups}, hep-th 9412059}.

The aim of this paper is to argue that all the possible quantum deformations
of $E(2)$ can be deduced from the analysis of Lie bialgebra structures on
$e(2)$. The paper is organized as follows.
In the chapter 2 we review obtained so far in the literature
quantum deformations of
$E(2)$ and $U(e(2))$.
In chapter 3 we present the classification of the Lie bialgebra
structures for $e(2)$. We obtain one one-parameter family of Lie bialgebra
structures and three separate "points".
Some
of them turn out to be coboundaries (i.e. they are determined by
classical $r$-matrices)
but
some are not of that kind. We show that quantum deformations described in
the chapter 2 give rise to all of them except one case. The missing
quantum deformation of $E(2)$ turns out to be a simplest one and will be
discussed in the chapter 4. In this chapter we
describe in detail how to derive
Lie-Poisson brackets
corresponding to Lie bialgebra structures which are not coboundaries.
In the new case of Lie Poisson brackets which did not yet appear in the
literature we obtain quantum group relations by changing Poisson brackets
into commutators. We also
calculate by duality the corresponding quantum deformation of
$U(e(2))$.
In chapter 5 we
conclude the paper
with some final remarks.

In our presentation we
will concentrate on algebra and coalgebra structures of
$E_q(2)$ and $U_q(e(2))$.
It is a trivial exercise to guess what is a form of antipode and
counit which make them Hopf algebras.

{\bf 2.} The first papers in the interesting us domain were
dedicated
to quantum deformations of enveloping algebra $U(e(2))$ \cgsta . These
deformations were obtained by applying the technique of contraction from
the standard deformation of $U(sl(2))$.
It turns out that there are two different quantum contractions
\ref\cgstb{E. Celeghini, R. Giachetti, R. Sorace and M. Tarlini, {\it
Contractions of quantum groups} in {\it Quantum Groups}, Lecture Notes
in Mathematics 1510, 221, (Springer Verlag, 1992)}:
\item{{\bf (A)}}
\eqn\relaa{[P_+, P_-] = 0,\qquad
[J, P_{\pm}] = \pm P_{\pm}}
\eqn\relac{\Delta (J) = J\otimes 1 + 1\otimes J}
\eqn\relad{\Delta (P_{\pm}) = P_{\pm}\otimes q^J + q^{-J}\otimes P_{\pm}}
\item{{\bf (B)}}
\eqn\relba{[P_1, P_2] = 0,\qquad
[J, P_1] = iP_2}
\eqn\relbc{[J, P_2] = - {i\kappa\over 2} sh {2P_1\over \kappa}}
\eqn\relbd{\Delta (J) = J\otimes e^{-P_1/\kappa} +
e^{P_1/\kappa}\otimes J}
\eqn\relbe{\Delta (P_1) = P_1\otimes 1 + 1\otimes P_1}
\eqn\relbf{\Delta (P_2) = P_2\otimes e^{-P_1/\kappa} +
e^{P_1/\kappa}\otimes P_2}

The deformation parameters are $q$ in the case ${\bf (A)}$ and $\kappa$ in
the case ${\bf (B)}$ and the classical limits are $q\rightarrow 1$ and $\kappa
\rightarrow\infty$.
It should
be perhaps mantioned that till now no general thery of contractions of quantum
groups exist. Some recent papers investigate such
contractions by analysing the Lie bialgebra level:
\ref\bhds{A. Ballesteros, F.J. Herranz, M.A. del Olmo and M. Santander,
J. Math. Phys. {\bf 36} (1995) 631} -
\ref\bghds{A. Ballesteros, N.A. Gromov, F.J. Herranz, M.A. del Olmo and
M. Santander, hep-th 9412083}.

Quantum $E_q(2)$ group has been discussed by many authors from
different
points of view
\cgsta -
\masa .
In order to fix the notation let us introduce the
following matrix representation of elements of $E(2)$:
\eqn\mate{g(c, a, b)=\pmatrix{\cos (c)&\sin (c)&a\cr
-\sin (c)&\cos (c)&b\cr 0&0&1}}
The matrix multiplication defines coproduct and
antipode
for $a, b$ and $c$. It turns out to be convenient to introduce the
complex notation:
\eqn\nota{\eta = a + ib,\qquad \bar\eta = a - ib\qquad {\rm and}\qquad e^{ic}}
The coproducts take the form:
\eqn\coma{\Delta (\eta ) = e^{-ic}\otimes \eta + \eta \otimes 1}
\eqn\comb{\Delta (\bar\eta ) = e^{ic}\otimes \bar\eta + \bar\eta
\otimes 1}
\eqn\comc{\Delta (e^{ic}) = e^{ic}\otimes e^{ic} }
There are many approaches to obtain quantum group relations for $E(2)$.
\wor - \bcgst .
They lead to two sets of relations:
\item{{\bf (A')}}
\eqn\graa{\eta \bar\eta = q^2 \bar\eta \eta}
\eqn\grab{\eta e^{ic} = q^2 e^{ic} \eta}
\eqn\grac{\bar\eta e^{ic} = q^2 e^{ic} \bar\eta}
or
\item{{\bf (B')}}
\eqn\grba{[e^{ic}, \eta ] = {1\over\kappa} (1 - e^{ic})}
\eqn\grbb{[e^{ic}, \bar\eta ] = {1\over\kappa} (e^{2ic} -
e^{ic})}
\eqn\grbc{[\eta, \bar\eta ] = {1\over\kappa} (\bar\eta + \eta )}
Coproducts for $\eta, \bar\eta$ and $e^{ic}$
are given in \coma -\comc .

It should be mentioned that the full Hopf
algebra duality has been demonstated for two discussed so far
deformations of the group $E(2)$ and $U(e(2))$. In the case {\bf (A)}
and {\bf (A')} it was shown in
\swz\
and in the case {\bf (B)} and {\bf (B')} in
\masa .

There exists still
another approach to the quantization of $E(2)$. The starting point is
the nonstandard (sometimes called Jordanian) quantum deformation of
$U(sl(2))$.
Following the general ideas it is possible to perform the contraction
from $U_q(sl(2))$ to $U_{\mu}(e(2))$
\ref\khpr{M. Khorrami, A. Shariati, M.R. Abolhassani and A.
Aghamohammadi, Mod. Phys. Lett. A {\bf 10} (1995) 873}.
The new deformation parameter is called $\mu$.
\item{{\bf (C)}}
\eqn\ncoa{[P_+,P_- ] = 0,\qquad
[J, P_+] = \mu sh {P_+\over \mu}}
\eqn\ncoc{[J, P_-] = - P_- ch {P_+\over \mu}}
\eqn\ncod{\Delta (P_+ ) = P_+ \otimes 1 + 1 \otimes P_+ }
\eqn\ncoe{\Delta (J ) = J \otimes e^{P_+/\mu} + e^{-P_+/\mu} \otimes J }
\eqn\ncof{\Delta (P_- ) = P_- \otimes e^{P_+/\mu} + e^{-P_+/\mu} \otimes
P_-}

The bad feature of this deformation of $U(e(2))$ is that it is strictly
speaking a deformation of the complex $e(2)$ algebra. This is seen in
the formulas \ncoa - \ncof . The operation $J^*=J$, $P_{\pm}^*=P_{\mp}$
is not a star operation in the Hopf algebra $U_{\mu}(e(2))$.

{\bf 3.} It is possible to give a complete
classification of Lie bialgebra structures for $e(2)$.
Let us introduce $e(2)$ Lie algebra with
generators
$P_1, P_2$ and $J$ satisfying  relations:
\eqn\etwo{[P_1, P_2] = 0,\qquad [J, P_1] = iP_2\qquad [J, P_2] = -iP_1.}

In the classification of Lie-bialgebra structures for $e(2)$ one
should take into account
its invariance under the following transformations: (i) $J$
$\rightarrow$ $J + \mu P_1 + \nu P_2$; (ii) $P_1$ $\rightarrow$ $\cos\beta
P_1 +
\sin\beta P_2$, $P_2$ $\rightarrow$ $-\sin\beta P_1 + \cos\beta P_2$;
(iii) $P_1$
$\rightarrow$ $\lambda P_1$, $P_2$ $\rightarrow$ $\lambda P_2$.
The complete
list of Lie bialgebra structures consist from:
\eqn\lba{\delta_1 (P_1) = s P_1\wedge J,\qquad
\delta_1 (P_2) = s P_2\wedge J,\qquad
\delta_1 (J) = 0}
\eqn\lbb{\delta_2 (J) = P_1\wedge P_2,\qquad
\delta_2 (P_1) = \delta_2 (P_2) = 0 }
\eqn\lbc{\delta_3 (P_2) = P_1\wedge P_2,\qquad \delta_3 (J) = P_1\wedge J,
\qquad \delta_3 (P_1) = 0}
\eqn\lbd{\delta_4 (P_1) = -i P_1\wedge P_2,\quad \delta_4 (P_2) =
P_1\wedge P_2,\quad
\delta_4 (J) = P_1\wedge J + i P_2\wedge J }
In the case \lba\ there is a one-parameter ($s$) family of Lie bialgebras.
Out of above four possibilities only the last two are coboundaries with
the classical $r$-matrices:
\eqn\baged{r_3 = J\wedge P_2}
\eqn\bagee{r_4 = J\wedge P_1 + i J\wedge P_2}
One could also write down more general form of
the classical $r$-matrices by adding terms
$\tau P_1\wedge P_2 +
\gamma (P_1\otimes P_1 + P_2\otimes P_2).$
This generalization will however turn out to be inessential if one
deduces the form of Lie-Poisson brackets out of $r$.
We find two Lie bialgebra structures
which are not a coboundary. It is interesting to stress that the
case of $D=2$ is very particular one. It was shown that for $D\geq 3$
all the Lie bialgebra structures of homogeneous groups built from
space-time rotations (with arbitrary signature) and translations
are coboundaries
\ref\zakc{S. Zakrzewski, {\it Poisson-Poincar\'e groups}, in {\it
Quantum Groups}, eds. J. Lukierski, Z. Popowicz and J. Sobczyk, PWN
Warszawa (Poland) 1995}.

It is easy to observe that
the deformation {\bf (A)} of chapter 2 corresponds
to $\delta_1$, {\bf (B)} corresponds to $\delta_3$ and {\bf (C)} to
$\delta_4$. On the other hand Lie bialgebra $\delta_2$ does not have yet
its quantum group counterpart.
\bigskip

{\bf 4.} When Lie bialgebra is not a coboundary the
computation of
Lie-Poisson brackets is not straightforward. The problem
is to solve
the cocycle equation for $\phi : G\rightarrow$ ${\cal G}$$\wedge$${\cal
G}$ where $G$ and ${\cal G}$ denote Lie group and its Lie algebra
\ref\drin{V.G. Drinfeld, Soviet. Math. Dokl. {\bf 27} (1983) 68}:
\eqn\coc{\phi (gh) = \phi (g) + g\phi (h) g^{-1}}
Let us discuss first the case of Lie bialgebra structure $\delta_1$.
The ``initial" conditions are ($\epsilon$ is an infinitesimal parameter)
\eqn\coina{\phi (1 + \epsilon P_1 + ...) = \epsilon s P_1\wedge J + ...}
\eqn\coinb{\phi (1 + \epsilon P_2 + ...) = \epsilon s P_2\wedge J + ...}
\eqn\coinc{\phi (1 + \epsilon J + ...) = 0}
It is understood that elements of $e(2)$ and $E(2)$ are given in
the $3$-dimensional representation and:
\eqn\rema{P_1=i\pmatrix{0&0&1\cr 0&0&0\cr 0&0&0},\qquad
P_2=i\pmatrix{0&0&0\cr 0
&0&1\cr 0&0&0},\qquad J=i\pmatrix{0&-1&0\cr 1&0&0\cr 0&0&0}}
The strategy is to find first $\phi$ on $1$-parameter subgroups
generated by $P_1, P_2$ and $J$. Let
\eqn\ans{\phi (e^{-iaP_1}) = A(a)P_1\wedge J + B(a) P_2\wedge J +
C(a) P_1\wedge P_2}
Cocycle equation \coc\ gives rise to the following set of
algebraic equations
\eqn\alea{2A(a) = A(2a)}
\eqn\alea{2B(a) = B(2a)}
\eqn\alea{2C(a) - aA(a) = C(2a)}
We assume further functions $A(a), B(a)$ and $C(a)$ to be analytic in
$a$. Taking into the account \coina\ one obtains
\eqn\xcoc{\phi (e^{-iaP_1}) = is \left(-a P_1\wedge J + {a^2\over 2}
P_1\wedge P_2\right)}
Using the same methods one calculates also:
\eqn\ycoc{\phi (e^{-ibP_2}) = is \left(-b P_2\wedge J + {b^2\over 2}
P_1\wedge P_2\right) }
\eqn\mcoc{\phi (e^{icJ}) = 0}
Group elements of $E(2)$ can be parameterized by
\eqn\grpar{g(a,b,c)=e^{-iaP_1}e^{-ibP_2}e^{icJ}=
\pmatrix{\cos c&\sin c&a\cr -\sin
c&\cos c&b\cr 0&0&1}}
for which using \coc\ one calculates
\eqn\ficoc{\phi (g(a,b,c)) = is \left( -aP_1\wedge J - bP_2\wedge J +
{a^2 + b^2\over 2}P_1\wedge P_2\right) }
Since one knows from the general theory that
\drin :
\eqn\bragen{\{f, g\} = \phi^{ab}\partial_a f\partial_b g}
the Lie-Poisson brackets for $a, b, c$ follow
\eqn\bbra{\{a, b\} = {is\over 2} (a^2 + b^2)}
\eqn\bbrb{\{a, \cos c\} = ias \sin c}
\eqn\bbrc{\{b, \cos c\} = ibs \sin c}
or, in the complex notation
\eqn\bnba{\{\eta , \bar\eta\} = s\eta\bar\eta}
\eqn\bnbb{\{\eta , e^{ic}\} = s\eta e^{ic}}
\eqn\bnbc{\{\bar\eta , e^{ic}\} = s\bar\eta e^{ic}}
These expressions should be compared with \graa -\grac .

Let us apply the same method to the second
non-coboundary
Lie bialgebra structure $\delta_2$:
\eqn\llb{\delta_2 (J)=P_1\wedge P_2,\qquad \delta_2 (P_1)=\delta_2 (P_2)=0}
Using once more the technique described above one obtains:
\eqn\lcoc{\phi (g(a,b,c)) = icP_1\wedge P_2}
After calculating Lie-Poisson brackets it turns out that the only
non-vanishing bracket is:
\eqn\lbra{\{\eta , \bar\eta \} = -2c}
One can check explicitely that it in fact satisfies the required condition
\eqn\check{\{\Delta (\eta ), \Delta  (\bar\eta )\} = \Delta (\{\eta ,
\bar\eta\}) = \Delta (-2c) = -2(c\otimes 1 + 1\otimes c)}
Naive quantization seems to applicable in this case so that one obtains
as the quantum relations
\eqn\lqua{[\eta , \bar\eta ] = ihc,\qquad [\eta , c] = [\bar\eta , c]=0}
where by $h$ we denoted the deformation (quantization) parameter. It is
instructive to find by duality the corresponding quantum deformation of
$U(e(2))$. After short computations one arrives at the following
structure
\item{{\bf (D)}}
\eqn\lqala{[J, P_1]= iP_2, \qquad [J, P_2] = -i P_1, \qquad [P_1, P_2]=0}
\eqn\lqalb{\Delta (P_1)=P_1\otimes 1 + 1\otimes P_1\qquad
\Delta (P_2)=P_2\otimes 1 + 1\otimes P_2}
\eqn\lqalc{\Delta (J) = J\otimes 1 + 1\otimes J + h(P_1\otimes P_2 -
P_2\otimes P_1)}

This is in fact the simplest possible quantum deformation of $U(e(2))$.
The antipode and counit are as in the undeformed case.
\bigskip

{\bf 5.}
We conclude that a theory of quantum deformations of the $D=2$
Euclidean group seems to be almost complete.
All the Lie bialgebra structures on
$e(2)$ can be quantized to Hopf algebras $U_q(e(2))$.
There is however still one interesting unsolved problem. It is
unknown whether for the quantum deformation \relba - \relbf\ universal
$R$-matrix exists. In many cases contraction
prescription can be applied to $R$-matrix giving rise to a finite
(usually after some manipulations) result
\ref\cgstc{E. Celeghini, R. Giachetti, R. Sorace and M. Tarlini, J.
Math. Phys. {\bf 32} (1991) 1155 and 1159} -
\ref\cgk{E. Celeghini, R. Giachetti and P. Kulish, J. Phys. {\bf A24}
(1991) 5675}. In the case {\bf (A)} discussed above the classical
$r$-matrix does not exist and the same must be true for the universal
$R$-matrix. In the case {\bf (B)} the situation is unclear.
Contraction of the universal $R$-matrix for
$U_q(sl(2))$ leads to divergent expressions.
Direct
computation shows however that, at least up to terms ${1\over
\kappa^{10}}$, an expression for $R$ can be found
\ref\ogiev{O. Ogievetsky, private communication}.
Such $R$ satisfies the condition $R\Delta (a) = \Delta '(a)R$ for all
the elements $a$ but it cannot satisfy Yang-Baxter equation as the
classical $r$-matrix satisfies only modified classical Yang-Baxter
equation.
\bigskip\noindent
ACKNOWLEDGEMENTS

I would like to thank dr. O. Ogievetsky for a helpful conversation.
\listrefs
\end

\ref\swz{P. Schupp, P. Watts and B. Zumino, Lett. Math. Phys. {\bf 24}
(1992) 141}

\ref\ohn{Ch. Ohn, Lett. Math. Phys. {\bf 28} (1992) 85}

\ref\wor{S.L. Woronowicz, Lett. Math. Phys. {\bf 23} (1991) 251; S.L.
Woronowicz, Comm. Math. Phys. {\bf 149} (1992) 637}

\ref\lnrt{J. Lukierski, A. Nowicki, H. Ruegg and V. Tolstoy, Phys. Lett.
B {\bf 264} (1991) 331}

\ref\frt{L.D. Faddeev, N.Yu. Reshetikhin and L.A. Takhtajan, Leningrad
Math. J. {\ bf 1} (1990) 193}

\ref\bcgsta{A. Ballesteros, E. Celeghini, R. Giachetti, R. Sorace and M.
Tarlini, J. Phys. A {\bf 26} (1993) 7495}

\ref\cgsta{E. Celeghini, R. Giachetti, R. Sorace and M. Tarlini, J.
Math. Phys. {\bf 31} (1990) 2548}

\ref\cgstb{E. Celeghini, R. Giachetti, R. Sorace and M. Tarlini, {\it
Contractions of quantum groups} in {\it Quantum Groups}, Lecture Notes
in Mathematics 1510, 221, (Springer Verlag, 1992)}

\ref\cgstc{E. Celeghini, R. Giachetti, R. Sorace and M. Tarlini, J.
Math. Phys. {\bf 32} (1991) 1155 and 1159}

\ref\drin{V.G. Drinfeld, Soviet. Math. Dokl. {\bf 27} (1983) 68}

\ref\zaka{S. Zakrzewski, J. Phys. {\bf A27} (1994) 2075}

\ref\zakb{S. Zakrzewski, Lett. Math. Phys. {\bf 22} (1991) 287}

\ref\zakc{S. Zakrzewski, {\it Poisson-Poincar\ e groups}, in {\it
Quantum Groups}, eds. J. Lukierski, Z. Popowicz and J. Sobczyk, PWN
Warszawa (Poland) 1995}

\ref\worb{S.L. Woronowicz, Rep. Math. Phys. {\bf 30 (2)} (1991) 259}

\ref\pw{P. Podle\ s and S.L. Woronowicz, {\it On the classification of
quantum Poincar\ e groups}, hep-th 9412059}

\ref\masa{P. Ma\ slanka, J. Math. Phys. \{ bf 35} (1994) 1976}

\ref\masb{P. Ma\ slanka, J. Phys. {\bf A27} (1994) 7099}

\ref\km{P. Kosi\ nski and P. Mas\ slanka,

\ref\zaua{P. Zaugg, J. Math. Phys. {\bf 36} (1995) 1547}

\ref\zaua{P. Zaugg, J. Phys. A {\bf 28} (1995) 2589}

\ref\bhds{A. Ballesteros, F.J. Herranz, M.A. del Olmo and M. Santander,
J. Math. Phys. {\bf 36} (1995) 631}

\ref\bghds{A. Ballesteros, N.A. Gromov, F.J. Herranz, M.A. del Olmo and
M. Santander, hep-th 9412083}

\ref\cgk{E. Celeghini, R. Giachetti and P. Kulish, J. Phys. {\bf A24}
(1991) 5675}

\ref\sob{J. Sobczyk, {\ $\kappa$ contraction from $SU_q$ to
$E_{\kappa}(2)$}

\ref\drber{V.G. Drinfeld, {\it Quantum groups}, Proc. ICM Berkeley 1986,
vol 1, 789}

\ref\lyu{V.V. Lyubashenko, Uspekhi Mat. Nauk, {\bf 41} (1986) 1587}

\ref\khpr{M. Khorrami, A. Shariati, M.R. Abolhassani and A.
Aghamohammadi, Mod. Phys. Lett. A {\bf 10} (1995) 873}

\ref\bchos{A. Ballesteros, E. Celeghini, F.J. Herranz, M.A. del Olmo and
M. Santander, J. Phys. A {\bf 28} (1995) 3129}

\ref\dvl{M. Dubois-Violette and G. Launer, Phys. Lett. {\bf B245} (1990)
175}

\ref\dmmz{E.E. Demichev, Yu.I. Manin, E.E. Mukhin and D.V. Zhdanowich,
Prog. Theor. Phys. Suppl. {\bf 102} (1990) 203}

\ref\aks{A. Aghamohammadi, M. Khorrami and A. Shariati, J. Phys. A {\bf
28} (1995) L225}

\ref\zaka{S. Zakrzewski, J. Phys. {\bf A27} (1994) 2075}.